\documentclass[11pt, a4paper, reqno]{amsart}
\usepackage[a4paper, margin=3.75cm]{geometry}
\usepackage{amssymb}
\usepackage{mathrsfs}
\usepackage{graphicx}

\usepackage[width=1.0\textwidth]{caption}

\usepackage[breaklinks]{hyperref}

\hypersetup{
unicode=true,
colorlinks=true,
linkcolor=blue,
citecolor=blue,
urlcolor=blue,
filecolor=blue,
bookmarksnumbered=true,
pdfstartview=FitH,
pdfhighlight=/N
}

\theoremstyle{plain}
\newtheorem{theorem}{Theorem}

\theoremstyle{definition}

\newtheorem{problem}{Problem}

\def\ZZ{\mathbb Z}
\def\Re{\operatorname{Re}}

\def\mdeg{\operatorname{deg}}

\def\supp{\operatorname{supp}}
\def\const{\mathrm{const}}

\def\mdeg{\operatorname{deg}}

\def\sZ{\mathscr Z}

\def\NN{\mathbb N}
\def\RR{\mathbb R}
\def\CC{\mathbb C}

\def\DD{\mathbb D}

\def\KK{\mathbf K}

\def\SZ{\mathscr Z}
\def\HH{\mathscr H}
\def\FF{\mathbf F}
\def\EE{\mathbf E}
\def\MM{\mathscr M}

\def\UU{\mathbb U}

\def\RS{\mathfrak R}

\def\zz{\mathbf z}
\def\mytt{\mathbf t}

\let\pfi\varphi
\let\geq\geqslant

\let\myh\widehat\let\myt\widetilde\let\myo\overline
\def\bad{\spaceskip=0.33emplus0.6emminus0.15em\immediate\write5{\string\bad}}

\def\dmu{\boldsymbol\mu}

\def\bGamma{\boldsymbol\Gamma}
\def\blambda{{\boldsymbol\lambda}}
\def\bmu{{\boldsymbol\mu}}
\def\myk{\mathfrak K}
\def\myD{\mathfrak D}

\begin{document}

\title{Structure of the Nuttall partition for some class of four-sheeted Riemann surfaces}

\author[Nikolay~R.~Ikonomov]{N.~R.~Ikonomov}
\address{Institute of Mathematics and Informatics, Bulgarian Academy of Sciences}
\email{nikonomov@math.bas.bg}
\author[Sergey~P.~Suetin]{S.~P.~Suetin}
\address{Steklov Mathematical Institute of Russian Academy of Sciences}
\email{suetin@mi-ras.ru}
\thanks{This research of the second author was carried out with the financial support of the Russian Foundation for Basic
Research (grant no.\ 18-01-00764).}

\date{March 8, 2021}

\begin{abstract}
The structure of a Nuttall partition   into sheets of some class of four-sheeted Riemann surfaces is
studied. The corresponding class of multivalued analytic functions is a~special class
of algebraic functions of fourth order generated
by the function inverse to the Zhukovskii function.
We show that in this class of four-sheeted Riemann surfaces, the boundary between the second and third sheets of the Nuttall partition
of the Riemann surface, is completely characterized in terms of an
extremal problem posed on the two-sheeted Riemann surface
of the function~$w$ defined by the equation $w^2=z^2-1$. In particular, we show that in this class of
functions the boundary between the second and third sheets does not intersect both the
boundary between the first and second sheets and the boundary between the third and fourth sheets.

Bibliography: \cite{Sue20} titles.

\end{abstract}

\maketitle

\markright{Structure of the Nuttall partition}

\setcounter{tocdepth}{1}
\tableofcontents


\section{Introduction and the statement of the main results}\label{s1}

\subsection{}\label{s1s1}
It is well known that a~global partition of a~Riemann surface of an algebraic function into ``sheets'' plays a~key role in the asymptotic
theory of Hermite--Pad\'e polynomials; see, first of all \cite{Nut84}, and also
\cite{Sta88},~\cite{KoPaSuCh17},~\cite{Kom20}.
Recently the  well-known Nuttall's conjecture (see \cite[\S\,3]{Nut84}), which has been open since 1984,
was proved in~\cite[\S~5, Lemma~5]{KoPaSuCh17}.
This conjecture claimed that, for the so-called Nuttall partition (with respect to some highlighted
point; see \cite{Nut84},~\cite{KoPaSuCh17},~\cite{Sue19}) of a~compact
$m$-sheeted Riemann surface into sheets, the complement of the ``topmost''
closed  $m$th sheet is always connected; that is, it is always a~domain on this Riemann surface. In~\cite{Sue18d},
this fact was used  as a~basis for the new approach to the implementation of the well-known Weierstrass
program (see~\cite{Ara03} and the bibliography given there) of effective  continuation of a~power series\footnote{In other words,
 here one speaks about the effective summation of a~power series beyond its convergence disk.
 By effective summation/continuation of a~power series one means its summation/extension
 via an (infinite) sequence of rational functions, in which each member of the sequence is
 constructed directly from a~finite number of coefficients of a~given  series.}
of a~given germ of a~multivalued analytic function. The verified Nuttall conjecture
naturally suggests the more detailed study of the topological properties of the Nuttall partition
of a~Riemann surface  into sheets. The need for this investigation stems, in particular, from the fact that the
further progress in the study of asymptotic properties
of the generalized Hermite--Pad\'e polynomials (which were introduced in~\cite{Kom20}, see also \cite{Sue18d})
requires a~more detailed knowledge of the topological structure of the Nuttall partition of a~compact Riemann surface into sheets.
However, many  properties  of a~Nuttall partition, which seem quite natural at first sight,
are still not rigorously justified in the general case. The study of such properties would allow one, in particular, to extend the Stahl theory of
the convergence of Pad\'e approximants, which holds, \textit{inter alia}, for infinitely-valued
analytic functions) to more general rational approximants of functions with the same properties
constructed on the basis of Hermite--Pad\'e polynomials.\footnote{Here and in what follows, we shall speak about Hermite--Pad\'e polynomials
 of type~I.}
 So, in parallel with the study of the asymptotic properties of Hermite--\allowbreak Pad\'e polynomials
 pertaining to algebraic functions and the corresponding compact  Riemann surfaces (see \cite{Nut84},~\cite{KoPaSuCh17},~\cite{Kom20}),
it is natural (in analogy with the Stahl theory) to investigate the asymptotic properties
of Hermite--Pad\'e polynomials for multivalued functions with finite number of singular
points on the Riemann sphere, but with the infinitely-sheeted Riemann surface.
In this direction, only the first steps have been made;  see \cite{RaSu13}, \cite{MaRaSu16},~\cite{Sue19}.

It is known  (see, in the first place, \cite{Sta97b}, and also \cite{Nut84} and~\cite{ApBuMaSu11}) that according to the Stahl theory,\footnote{By the Stahl theory
one usually means the series of results obtained by Stahl  in 1985--1986s on the convergence
of Pad\'e approximants to multivalued analytic functions;
see~\cite{Sta97b}.}
with each germ
$f_{z_0}$ (considered at a~point $z_0\in\myh{\CC}$) of a~multivalued analytic function~$f$ with a~finite number of singular
points on the Riemann sphere $\myh{\CC}$ one can uniquely associate a~hyperelliptic Riemann surface, whose first sheet\footnote{We assume,  unless otherwise stated, that all the sheets
of a~Riemann surface are open. Note that it is well known that, for a~Nuttall partition
of a~Riemann surface  into sheets, these sheets may fail to be connected.} corresponds to the
uniquely defined (from the given germ  $f_{z_0}$) Stahl domain, which is
the maximal domain of convergence\footnote{The convergence of Pad\'e approximants is usually understood
in the sense of the convergence in logarithmic capacity on compact  subsets of the Stahl domain.} of the diagonal Pad\'e approximants.

In the new approach of~\cite{Sue18d} to the solution of the
problem of effective analytic continuation of a~given multivalued
function, when the Hermite--Pad\'e polynomials  are considered in leu of
the Pad\'e polynomials, one should deal, instead of a~single (multivalued) function~$f$ (or, in
a~different terminology, a~family of functions $[1,f]$), with the pair of functions
$f,f^2$ (or, in other words, with the family of functions $[1,f,f^2]$). According to~\cite{Sue18d} (see also \cite{Nut84} and~\cite{KoPaSuCh17}),
this approach in principle allows one to double   (in comparison
with the Stahl domain) the domain in which the values of the function~$f$
can be efficiently (that is, in terms of rational  functions) recovered from a~given germ.
The corresponding generalizations within the framework of this new approach are
given in~\cite{Sue18d} also for the family $[1,f,f^2,f^3]$ (see also \cite{Sue18b},~\cite{Sue20}).
In this case, the domain of effective recovery of an analytic function increases already by three times
in comparison with the Stahl domain. Namely, the corresponding Nuttall domain
lies on the Riemann surface, and hence in a~certain sense this domain
is a~three-sheeted covering of the Stahl domain, which lies on the Riemann sphere.
A~much more general case of the family $[1,f,\dots,f^m]$ ($m\in\NN$ is arbitrary,
$f$~is an algebraic  function of order $(m+1)$) was considered by Komlov~\cite{Kom20}.

All that was said above, was discussed in \cite{Sue18d} on an example
of multivalued functions $f$  from the class introduced earlier by the author of the present paper (see~\cite{Sue18})\,---\,this being the
class of functions of the form
\begin{equation} f(z):=\left[\left(
A-\frac1{\pfi(z)}\right)\left(B-\frac1{\pfi(z)}\right)\right]^{-1/2}, \quad z\in
D:=\myh{\CC}\setminus{E}, \quad E:=[-1,1],
\label{1}
\end{equation}
where
$\pfi(z):=z+(z^2-1)^{1/2}$ is the function inverse to the Zhukovskii function, $1<A<B<\infty$;
here and what follows we choose the branch of the function  $(\,{\cdot}\,)^{1/2}$ such that
$(z^2-1)^{1/2}/z\to1$ as $z\to\infty$ (for details, see \cite{Sue18} and
\S\,\ref{s1s2} below).

However, for an arbitrary multivalued function (in particular,
for an algebraic functions of order  $\geq4$), the question of asymptotic
properties of the Hermite--\allowbreak Pad\'e polynomials for the family $[1,f,f^2]$
remains open in the general case.\footnote{A~corresponding general result would give a~natural
extension of Stahl's  theorem to the case of Hermite--Pad\'e polynomials  for the family
$[1,f,f^2]$ (in place of the family  $[1.f]$).} The corresponding results
in this direction have been obtained to date only in very special cases;
see, for example,~\cite{RaSu13},~\cite{LoLo18},~\cite{Sor20}).
One of the main reasons is as follows. In the Stahl theory, the
domain of convergence of Pad\'e approximants  is completely characterized\footnote{Up to a~set of zero capacity.}
by the property that the boundary of this domain has the smallest capacity in the class of all the domains
admissible\footnote{We recall that a domain $G\ni z_0$ on the Riemann sphere
is called admissible for a~given germ $f_{z_0}$ if this germ extends from the point~ $z_0$ to this domain as a~meromorphic (single-valued
analytic) function.}
for a~given germ of a~multivalued function. So, the maximal Stahl domain corresponds to the solution of a~certain 
extremal problem (the problem of an admissible compact set of minimal capacity).
For the family $[1,f,f^2]$ (or, in a~different terminology, for the pair of functions $f,f^2$), where
the function~$f$ lies in the same class,\footnote{Of course, the case of a~hyperelliptic function is excluded from
consideration.} which was considered by Stahl,
all the attempts to pose the corresponding extremal problem in a~sufficiently general case
have proved futile.

We can mention several extreme problems associated with the study of the asymptotic behavior of
Hermite--Pad\'e polynomials, but which are
suitable only for certain special cases and extending only slightly the Stahl theory to the setting of Hermite--Pad\'e polynomials
(see \cite{Apt08},~\cite{RaSu13},~\cite{ApTu16},~\cite{Sue18d},~\cite{Sue19},~\cite{Sor20}).
Of course, the case with four functions $[1,f,f^2,f^3]$  is even more involved.
In this setting, some or other characterization of the Nuttall partition
of a~Riemann surface into sheets in terms of an appropriate extremal problem
(even in sufficiently simple cases when the structure of the partition of a~Riemann surface  into sheets
is fairly simple and quite  natural) would be of great value from the point of view of the further
development of the Stahl theory and its extension to the Hermite--Pad\'e polynomials setting.

This is precisely the purpose of the present paper\,---\,for a~function of the form \eqref{2}
to characterize the Nuttall partition of the corresponding four-sheeted Riemann surface in terms of the
extremal problem posed on the two-sheeted Riemann surface of the function  $w^2=z^2-1$
(more precisely, here we speak about the characterization of  the boundary between the second and third sheets
of the Nuttall partition). Namely, in the present paper
we consider some class of multivalued analytic functions which have only second-order branch points.
Moreover, it is assumed that geometrically the branch points are symmetric with respect to  the real line.
In other words, 
all $\alpha_j\in\{1/2,-1/2\}$ in representation \eqref{2}. In the present paper, this class will be denoted by~$\sZ$.
In accordance with the general result of \cite{KoPaSuCh17}, under the Nuttall partition (with respect to the
 the point at infinity $\zz=\infty^{(1)}$) of the corresponding Riemann surface into sheets, the
 complement of the fourth (last, the ``topmost'') sheet is always connected.  Furthermore, from
Komlov's general results \cite{Kom20} it follows that
using the Hermite--Pad\'e polynomials (of type~I and~II, respectively, and the generalized ones; see~\cite{Sue18d}
and~\cite{Kom20}) the values of a~function $f\in\sZ$ can be recovered on the first, second, and third sheets
of the Riemann surface $\RS_4(f)$. In the present paper, under certain geometric constraints on the
arrangement of branching points of a~function~$f$ from the class~$\sZ$, we give a~characterization
of the  boundary between the second and third sheets of its Riemann surface $\RS_4(f)$ in terms of
a~certain extremal problem posed now not on the Riemann sphere~$\myh{\CC}$, but on the
two-sheeted Riemann surface $\RS_2(w)$ of the function~$w$, which is given by the equation $w^2=z^2-1$.  Note that this extremal problem
calls for the
characterization of an ``admissible'' compact set located on the Riemann surface $\RS_2(w)$ and
having the property that the Green function of its complement (assuming that this complement is connected)
maximizes the Robin constant at one of the two points at infinity  located on the Riemann surface $\RS_2(w)$.
Thus, this extremal problem is quite analogous to the classical problem of an admissible
compact set of minimal capacity, which was posed (and solved) by Stahl  on the Riemann sphere.

In the last \S\,\ref{s4}, we compare the extremal Problem~\ref{pro1} considered here  and
the extremal problem of~\cite{Sue18b} on the existence of a~three-sheeted Nuttall-partitioned Riemann surface
associated with a~given multivalued analytic function.

Thus, in the present paper, we further extend the new approach (which was proposed by the
author in~\cite{Sue18b}) to the study of asymptotic properties of Hermite--\allowbreak Pad\'e polynomials  for multivalued functions.
This approach is based on the extremal equilibrium problem posed not on the Riemann sphere, but
instead on the Riemann surface $\RS_2(w)$ (further advances in this problem were made in ~\cite{Sue19b},~\cite{IkSu20}).

To conclude our introduction, we comment on the condition to the effect that all the exponents in \eqref{2}
are  $\pm1/2$. It is well known that the first results (obtained prior to the completion of the general Stahl theory)
on the convergence of Pad\'e approximants in the class
of multivalued analytic functions were obtained by Nuttall for the class of hyperelliptic functions and
their natural generalizations; that is, for functions
having only second-order branch points;
see \cite{NuSi77},~\cite{Nut84}, and
also~\cite{Sue00},~\cite{ApBuMaSu11},~\cite{ApYa15}.

It worth pointing out again that the problem of effective recovery of a~function
$f\in\sZ$ on three Nuttall sheets of its Riemann surface $\RS_4(f)$, which was considered in~\cite{Sue18d},
is a~very particular case of the general problem of the recovery of the values of an algebraic function~$f$
of order $m+1$ on the first $m$-Nuttall  sheets of its Riemann surface  $\RS_{m+1}(f)$ (this problem was considered by
Komlov in~\cite{Kom20}). From the properties
of the Nuttall partition of the Riemann surface $\RS_4(f)$ into sheets, which we obtain below, it follows that
the results of~\cite{Kom20} on the convergence of rational approximants constructed on the basis of Hermite--Pad\'e polynomials
are also valid in the case considered here.

\subsection{}\label{s1s2}

We will now give the necessary definitions and notation.

Let, as above, $\pfi(z):=z+(z^2-1)^{1/2}$  be the function inverse to the Zhukovskii function, which is
meromorphic and single-valued in the domain $D:=\myh{\CC}\setminus{E}$,
${E}:=[-1,1]$, and which maps the domain~$D$ onto the exterior  $\UU:=
\{\zeta\in\myh{\CC}:|\zeta|>1\}$ of the unit disk
$\DD:=\{\zeta\in\CC:|\zeta|<1\}$. Here and in what follows, we choose the  branch
of the function $(\,{\cdot}\,)^{1/2}$ such that $(z^2-1)^{1/2}/z\to1$ as $z\to\infty$.

Consider the class $\sZ$ of functions $f$ defined for $z\in D$ by the explicit representation
\begin{equation}
f(z)=\prod_{j=1}^{2p}\left(A_j-\frac1{\pfi(z)}\right)^{\alpha_j},\quad z\in D,
\label{2} \end{equation}
where $p\in\NN$ is a~natural number, all the exponents  $\alpha_j\in\{1/2,-1/2\}$,
$\sum_{j=1}^{2p}\alpha_j=0\pmod{\ZZ}$, all the quantities  $A_j$ are pairwise different,
$|A_j|>1$, and the set $\{A_1,\dots,A_{2p}\}$ is symmetric
with respect to the real line.  The function $f$ of the form~\eqref{2} has the following properties:

1) to the function~$f(z)$, $z\in D$, there corresponds the germ $f_\infty$ holomorphic at the
point $z=\infty$;

2) $f$ is an algebraic function of fourth order;

3) the function $f(z)$, $z\in D$, extends analytically  along any
path lying on the Riemann sphere $\myh{\CC}$ and not intersecting the set of points
$\Sigma=\Sigma(f)=\{\pm1,a_j,j=1,\dots,2p\}$, where $a_j=(A_j+1/A_j)/2\notin{E}$.

4) each point of the set~$\Sigma$ is a~second-order branch point of the function~$f$.

By $\RS_4(f)$ we denote the four-sheeted Riemann surface of the function  $f\in\sZ$; the Riemann surface
$\RS_4(f)$ has genus~$p-1$.

Note that the class of functions of the form~\eqref{2} is a~subclass of the more general class
studied in~\cite{Sue18d} and~\cite{Sue18}.

Consider the two-sheeted Riemann surface $\RS_2(w)$ of an algebraic function $w$ defined by
the equation $w^2=z ^2-1$. A~point $\zz\in\RS_2(w)$ on this Riemann surface  is the pair
$\zz=(z,w)$. Let $\pi_2$ be the corresponding canonical  projection of $\RS_2(w)$
onto the Riemann sphere $\myh{\CC}$: $\pi_2\colon\RS_2(w)\to\myh{\CC}$,
$\pi_2(\zz)=z$. We set $\bGamma:=\{\zz\in\RS_2(w):\pi_2(\zz)\in{E}\}$
(that is, $\bGamma=\pi^{-1}_2({E}) $). Let
$\infty^{(1)}\in\pi^{-1}_2(\infty)$ be a~point of the set
$\pi^{-1}_2(\infty)$ such that $w/z\to1$ as $\zz\to\infty^{(1)}$;
let $\infty^{(2)}$ be the second point of the set  $\pi^{-1}_2(\infty)$ ($w/z\to-1$
as $\zz\to\infty^{(2)}$). Let $\RS_2^{(1)}$ and $\RS_2^{(2)}$ be, respectively, the first and second  (open) sheets of the Riemann surface $\RS_2(w)$,
$z^{(1)}=(z,(z^2-1)^{1/2})$ and let $z^{(2)}=(z,-(z^2-1)^{1/2})$ be points lying, respectively,  on the first and  second sheets.

Let $\Phi(\zz):=z+w$ and $g(\zz,\infty^{(1)},\infty^{(2)}):=\log|z+w|
=\log|z\pm(z^2-1)^{1/2}|=\log|\Phi(\zz)|$, $\zz\in\RS_2(w)$,  be the bipolar Green function for the Riemann surface  $\RS_2(w)$
(see \cite{Chi18},~\cite{Chi19}). Note that
$g(\zz,\infty^{(1)}, \infty^{(2)})\equiv0$ for $\zz\in\bGamma$.

So,
$g(z^{(1)},\infty^{(1)},\infty^{(2)})=g_E(z,\infty)$ is the Green function for the domain
 $D=\myh{\CC}\setminus{E}$ (note that ${E}$ is a~compact Stahl set for
$f\in\sZ$, $D$~is the corresponding Stahl domain). For
$\zz=z^{(2)}$ we have
$\log|z+w|=\log|z-(z^2-1)^{1/2}|=-\log|z+(z^2-1)^{1/2}|$, and hence the function
$\eta(\zz):=-g(\zz,\infty^{(1)},\infty^{(2)})$ defines a~Nuttall partition
of the Riemann surface $\RS_2(w)$ into sheets; see \cite[\S~3]{Nut84},~\cite{KoPaSuCh17}. Namely, $\eta(z^{(1)})<\eta(z^{(2)})$ for $z\in D$.

Let $f_\infty$ be a germ of a~function~$f\in\sZ$ defined in the domain~$D$ under the above conditions.
The germ $f_\infty$ is lifted to the point  $\zz=\infty^{(1)}$ lying on the first sheet $\RS^{(1)}_2$
of the Riemann surface $\RS_2(w)$. The corresponding germ $f_{\infty^{(1)}}$ extends from the point $\zz=\infty^{(1)}$
to the entire first sheet of the Riemann surface $\RS_2(w)$. From the first sheet to the second sheet this germ
extends as a~multivalued function. So, in order to define a~single-valued meromorphic
extension of $f_{\infty^{(1)}}$ on the second sheet, one needs to introduce an appropriate family of
admissible compact sets (``cuts'') $\KK$ lying on the second sheet,
$\KK\subset\RS_2(w)$. Namely, we consider compact sets such that:

1) the compact set $\KK=\KK^{(2)}\subset\RS_2^{(2)}(w)$;

2) the compact set $\KK$ does not separate the Riemann surface $\RS_2(w)$; that is, the complement of $\KK$ is connected,
$\RS_2(w)\setminus\KK=D(\KK)$  is a~domain on $\RS_2(w)$, and $\infty^{(1)}\in
D(\KK)$;

3) the germ $f_{\infty^{(1)}}=f_\infty$ extends from the point
$\zz=\infty^{(1)}$ in the domain $D(\KK)$ as a~meromorphic (single-valued
analytic) function $f\in\MM(D(\KK))$.

For a given germ $f_\infty$ of the function~$f\in\sZ$, the class of all compact sets
$\KK\subset\RS_2(w)$ satisfying the above conditions~1)--3) will be denoted by
$\myk(f_\infty)$. Compact sets from this class will be called {\it
admissible}.  It is clear that the family $\myk(f_\infty)$ is nonempty.

It is also clear that the compact set $\KK\in\myk(f_\infty)$ is a~nonpolar set
(for this concept on a~Riemann surface, see~\cite{Chi19}). So, the Green function $g_{\KK}(\zz,\infty^{(1)})$ for the domain  $D(\KK)$
with singularity at the
point $\zz=\infty^{(1)}$ is well defined; see \cite{Chi19},~\cite{Chi19b},~\cite{Chi20}. At the point
$\zz=\infty^{(1)}$ we introduce the local coordinate~$\xi$. For an arbitrary
compact set $\KK\in\myk(f_\infty)$, we have (in terms of this coordinate)
\begin{equation}
g_{\KK}(\zz,\infty^{(1)})=\log\frac1{|\xi|}+\gamma(\KK)+o(1),
\quad\text{as}\quad\zz\to\infty^{(1)}.
\label{3}
\end{equation}
The constant $\gamma(\KK)$ in~\eqref{3}, which depends on the local
coordinate~$\zeta$, is called the Robin constant (with respect to this coordinate) at the point $\zz=\infty^{(1)}$.

In the class of compact sets $\myk(f_\infty)$, we pose the following extremal problem.

\begin{problem}\label{pro1} Prove the existence and characterize
an admissible compact set $\FF\in\myk(f_\infty)$ such that
\begin{equation}
\gamma(\FF)=\sup_{\KK\in\myk(f_\infty)}\gamma(\KK).
\label{4}
 \end{equation}
\end{problem}

It is clear that even though the quantity $\gamma(\FF)$ depends on the local
coordinate~$\xi$, the extremal  compact set (provided it exists) does not depend on the local coordinate.

The following result holds.

\begin{theorem}\label{the1}
There exists a~unique compact set $\FF\in\myk(f_\infty)$ satisfying condition~\eqref{4},  that is,
$\gamma(\FF)=\max\limits_{\KK\in\myk(f_\infty)}\gamma(\KK)$.
The compact set~$\FF$ consists of a~finite number of analytic arcs and
is completely characterized by the following $S$-property\footnote{For the definition of the
$S$-property, see \cite{GoRa87},~\cite{GoRaSu11},~\cite{Rak12},~\cite{Rak18}
and there references given there.}
\begin{equation} \frac{\partial
g_{\FF}(\zz,\infty^{(1)})}{\partial n^{+}} =\frac{\partial
g_{\FF}(\zz,\infty^{(1)})}{\partial n^{-}}, \qquad \zz\in\FF^{\circ},
\label{6}
\end{equation}
where $\FF^{\circ}$ is the union of all open arcs whose closures comprise the compact set $\FF$,
and $\partial/\partial
n^{\pm}$ are the normal derivatives at a point~ $\zz\in\FF^{\circ}$ from the opposite sides of~$\FF^{\circ}$.
 \end{theorem}


Note that the case  $\FF\ni\infty^{(2)}$ is not excluded from consideration.

In what follows, in the proof of Theorem~\ref{the2} we shall assume that the following
``general position'' condition  is satisfied: the compact set~$\FF$ consists precisely of~$p$
pairwise disjoint analytic arcs, that is, the compact set~$\FF$ contains no Chebotarev points.
This assumption is quite analogous to that adopted in the first papers on
asymptotic properties of Pad\'e polynomials;
see \cite{NuSi77},~\cite{Nut84},~\cite{Sue00}. Note that this assumption is equivalent to
saying that all zeros of the polynomial $V_{m-2}$, which appear in
a~characterization of the Stahl compact set  in terms of a~variational method, are of even multiplicity; see~\eqref{9}.


Theorem~\ref{the1} follows in part from the  general results of the book \cite[Ch.~8, \S~4]{ScSp54}.
However in~\cite{ScSp54} the proof of a~more general result is based on the Schiffer variational method.
Nevertheless, below we shall give our independent
proof of Theorem~\ref{the1}, because from this proof it will be  possible to derive an additional information on the
structure of the compact set~$\FF$ and its relation to some quadratic differential. Namely, it will be shown that the
compact set~$\FF$ consists of the critical trajectories of this quadratic differential.
This information will be required below in the proof of Theorem~\ref{the2}.
Theorem~\ref{the1} will be proved in~\S~\ref{s2}.

Note that Problem \ref{pro1} is quite similar to the classical problem on an admissible compact set of minimal  capacity
from the Stahl theory of Pad\'e approximation. In the classical case,
for an extremal compact set $S$
on the Riemann sphere, the quantity $e^{-\gamma(S)}$ is the capacity of this compact set  (with respect to
the point at infinity $z=\infty$).
It is well known  (see \cite{Sta97b}, and also~\cite{ApBuMaSu11}) that the Stahl compact set~$S$
satisfies the characteristic relation
\begin{equation}
\frac{g_S(z,\infty)}{\partial n^{+}}= \frac{g_S(z,\infty)}{\partial
n^{-}},\qquad z\in S^{\circ}, \label{7}
\end{equation}
which is quite similar to \eqref{6}. Here,  $S^{\circ}$ is the union of open arcs whose closures comprise
the Stahl compact set~$S$ and
$g_S(z,\infty)$ is the Green function for the domain  $\myh{\CC}\setminus{S}$.
We recall that the Stahl compact set~$S$ (an admissible  compact set of minimal capacity)
consists of the closures of critical trajectories of the quadratic differential. Namely,
the following relations hold:
\begin{align}
S&=\biggl\{z\in\CC:-\frac{V_{m-2}(z)}{B_m(z)}\,dz^2>0\biggr\}, \label{8}\\
g_S(z,\infty)&=\Re\int_{b_0}^\infty
\sqrt{\frac{V_{m-2}(\zeta)}{B_m(\zeta)}}\,d\zeta, \qquad z\in
\myh{\CC}\setminus S; \label{9}
\end{align}
here $B_m(z):=\prod\limits_{b\in\Sigma}(z-b)$, $\Sigma=\Sigma(f)$  is the set of singular  points of a~multivalued function
$f\in\HH(\infty)$, $\#\Sigma=m<\infty$, $\mdeg V_{m-2}=m-2$,
$V_{m-2}(z)=z^{m-2}+\dotsb=\prod\limits_{j=1}^{m-2}(z-v_j)$, $v_j$ are the Chebotarev  points
of the compact set~$S$.

So, with a~given germ $f_\infty$ of a~multivalued analytic function
$f$ with finite number of branching points on the Riemann sphere  one associates in a~unique  way
a~(unique) two-sheeted hyperelliptic Riemann surface
$\RS_2(f_\infty)$ defined by the equation $w^2=V_{m-2}(z)/B_m(z)$. Since by the above
the rational  function\footnote{Note that the polynomials
$V_{m-2}$ and $B_m$ may have common zeros. Hence
$V_{m-2}/B_m=V^{*}_{m-2}/B^{*}_m$, where now the polynomials $V^{*}_{m-2}$ and $B_{m}^*$
are relatively prime polynomials. Moreover, the polynomial $V^{*}_{m-2}$ may have zeros
of even multiplicity. As a~result, the two-sheeted Riemann surface $\RS_2(f_\infty)$ is defined in fact
by the quadratic equation $w^2=\myt{V}^{*}_{m-2}B^{*}_m$, where the polynomial
$\myt{V}^{*}_{m-2}$ is obtained from the polynomial  $V^{*}_{m-2}$ by removing the zeros of even multiplicity.}
$V_{m-2}(z)/B_m(z)$ is uniquely defined from the germ $f_\infty$, the function $w$ is also uniquely defined from the original germ.
This hyperelliptic Riemann surface is known as the {\it
$($Stahl$)$ associated} surface with the germ$f_\infty$. With the exception of the case when the original germ $f_\infty$
is a~germ of a~hyperelliptic  function, $f_\infty$~extends not on the entire associated Riemann surface, but also to the
first sheet of this surface. For a~further extension of this
germ~$f_\infty$ as a~single-valued analytic function on the second sheet of this Riemann surface, we need to organize the
corresponding cuts (in fact, the family of such cuts forms an admissible
compact set, and the corresponding family of compact sets forms the family of admissible compact sets).
Nevertheless, it turns out the so-called
{\it strong} asymptotics of the Pad\'e polynomials is characterized precisely in terms pertaining to this
two-sheeted Riemann surface $\RS_2(f_\infty)$ (see \cite{Nut84},~\cite{Sue00},~\cite{ApYa15}).
So, our approach, in which properties of extremal compact sets pertaining to Hermite--Pad\'e polynomials
are studied  with the help of the results obtained earlier in the Stahl theory and
its further advances made by Stahl himself and other researchers
(see \cite{Sue00},~\cite{Sue18b},~\cite{Sue19},~\cite{IkSu20}), is also quite natural.
The new approach proposed in \cite{Sue18b} has proved instrumental in
delivering, for the class of functions
of the form~\eqref{2}, some new and previously available results related to
the Hermite--Pad\'e polynomials in terms of the scalar equilibrium problem (posed on a~two-sheeted Riemann surface),
rather than in terms of the generally accepted equilibrium problem (posed on the Riemann sphere).
This scalar approach was further advanced in~\cite{Sue18},~\cite{Sue18b},~\cite{Sue19},~\cite{IkSu20} for a~pair of
functions forming a~Nikishin system  (cf.~\cite{RaSu13}). In particular, this also pertains to the pair
of functions $f,f^2$, which, as was shown by the author of the present paper in~\cite{Sue18},
forms a~Nikishin system (under a~minimal extension of this classical concept).
With this approach, the extremal problems on the corresponding two-sheeted Riemann surface were posed and solved
(see also \S\,\ref{4} below). This approach leads naturally to the (now three-sheeted) Riemann surface
which is Nuttall-associated with the original germ~$f_\infty$. Namely, this germ $f_\infty$,
as lifted to the point  $\zz=\infty^{(1)}$  as a~germ $f_{\infty^{(1)}}$, extends as a~(single-valued) meromorphic function
from this point to  the domain defined as the complement\footnote{The fact that this complement is a~domain
is a~part of the well-known Nuttall  conjecture~\cite{Nut84} of 1984.
This conjecture was proved in 2017 in~\cite{KoPaSuCh17}.} of the closure of the ``topmost'' third sheet.
However, for further single-valued extension of the germ  $f_{\infty^{(1)}}$ to the third sheet of this
Riemann surface, corresponding cuts are required  (cf.~\cite{RaSu13} and~\cite{Sta88}).

In the present paper, we generally adhere to the scalar approach, which, however, is developed in
a~sightly different situation. Namely, we formulate and solve the
extremal problem now for the family of four functions $[1,f,f^2,f^3]$, where
$f\in\sZ$.  Note that if $p=1$ in \eqref{2}  and if $1<A<B$, then according to~\cite{Sue18} the three functions
$f,f^2,f^3$ form a~Nikishin system.

For functions of the form~\eqref{2}, the Stahl compact set is the closed interval  ${E};=[-1,1]$.
Likewise, the two-sheeted Riemann surface Stahl-associated with an arbitrary
function $f\in\sZ$ is the Riemann surface $\RS_2(w)$ of the function~$w$ defined by the
equation $w^2=z^2-1$. On this Riemann surface $\RS_2(w)$, the extremal Problem~\ref{pro1} is formulated
 (cf.~\cite{RaSu13},~\cite{IkSu20}).

\subsection{}\label{s3s1}

Let $\RS_4(f)$ be the four-sheeted Riemann surface of a~function  $f\in\sZ$ and let
$\pi_4\colon\RS_4(f)\to\myh{\CC}$ be the corresponding canonical  projection
(see Fig.~\ref{fig_1} in the case $p=1$ and $1<A<B$ in \eqref{2}).

Let $f_\infty\in\HH(\infty)$ be the above germ of a~function~$f\in\sZ$ (that is,  $f_\infty$ extends
holomorphically  from the point $z=\infty$
to the Stahl domain  $D=\myh{\CC}\setminus{E}$). We shall assume that  the first
(open) sheet $\RS_4^{(1)}$ of the Riemann surface $\RS_4(f)$ is chosen so that
$\RS_4^{(1)}\simeq D$ and the mapping $\pi_4\colon\RS_4^{(1)}\to D$
is biholomorphic. Let $\infty^{(1)}:=\pi_4^{-1}(\infty)\cap\RS_4^{(1)}$.
Then the germ $f_\infty$ is lifted to the point  $\infty^{(1)}\in\RS_4^{(1)}$ and
extends everywhere to the Riemann surface $\RS_4$ as a~single-valued  meromorphic function.
In what follows, we will identify the germs $f_\infty$ and $f_{\infty^{(1)}}$, retaining the above notation  $f_\infty$.

Let us now define the global partition of the Riemann surface $\RS_4(f)$ into sheets as follows.
The Green function $g_{\FF}(\zz,\infty^{(1)})$ for the domain
$\myD_1:=\RS_2(w)\setminus\FF$ (the complement of the extremal compact set
$\FF=\mathbf F^{(2)}\subset\RS_2(w)$) has symmetric boundary
\eqref{6}, and moreover, by the above assumption, the boundary consists of a~finite number of disjoint analytic arcs.
Besides, by the definition of the Green function, we have
$g_{\FF}(\zz,\infty^{(1)})\equiv0$ for $\zz\in\FF$.  By arranging cuts on
$\RS_2(w)$ we may consider these arcs as two-sided arcs  (with the exception of the  end-points $\FF\setminus\FF^{\circ}$).
From the surface
$\mathfrak P_1\simeq\RS_2(w)\setminus\FF$, we construct a~four-sheeted Riemann surface
$\RS_4=\RS_4(f_\infty)$ as follows. Consider the second copy
$\mathfrak P_2$ of this surface ($\mathfrak
P_2\simeq\RS_2(w)\setminus\FF$)  and ``glue'' it together with the first copy by identifying the opposite sides (``edges'')
of the new cuts on these two copies of the surface $\RS_2(w)\setminus\FF$. We have
$\mathfrak P_1\simeq\mathfrak P_2$, and hence the points $\zz_1\in\mathfrak
P_1$ and $\zz_2\in\mathfrak P_2$ are in a~one-to-one correspondence, which we denote by the sign~``$\simeq$'':
$\zz_1\simeq\zz_2$. Note that  $\pi_2(\zz_1)=\pi_2(\zz_2)$.

In accordance with the above ``general position'' assumption to the effect that the arcs comprising the
compact set~$\FF$ are disjoint, it is clear that the four-sheeted Riemann surface $\RS_4(f_\infty)$ thus obtained coincides
with the Riemann surface $\RS_4(f)$. Thus, the original germ $f_\infty$ extends to the entire
surface $\RS_4$ as a~single-valued meromorphic function.  We also note that the fact
the four-sheeted Riemann surface $\RS_4$ thus constructed coincides with the Riemann surface
$\RS_4(f)$ of the function~$f\in\sZ$ can be also derived directly from the proof of Theorem~\ref{the1}.

Let us now define the global partition of the Riemann surface  $\RS_4(f)$ into sheets as follows.  In
accordance with the definition of the Green function, we have
$g_{\FF}(\zz,\infty^{(1)})=0$ for $\zz\in\FF$. From this equality and the symmetry condition
\eqref{6} of the $S$-compact set~$\FF$ it follows that the Green function
$g_{\FF}(\zz,\infty^{(1)})$ extends to a~harmonic function on the second part of the Riemann surface
$\RS_4(f)$ with the same values but with different sign
(``minus'' instead of ``plus''). Namely,
\begin{equation}
g_{\FF}(\zz_2,\infty^{(1)})=-g_{\FF}(\zz_1,\infty^{(1)})\quad\text{for}\quad\zz_1\in\mathfrak P_1,\quad\zz_2\in\mathfrak P_2.
\label{10}
\end{equation}
Thus, the compact set $\FF$, provided that it is composed of precisely $p$~disjoint analytic arcs
and these arcs are considered as two-sided cuts (with two ``edges') on the Riemann surface
$\RS_2(w)$, is the boundary between the second and third sheets of the Riemann surface
$\RS_4(f)$. The boundary between the first and second sheets goes along the curve
$\bGamma^{(1,2)}$ such that $\pi_4(\bGamma^{(1,2)})={E}$.  The boundary
between the third and fourth sheets goes along the curve  $\bGamma^{(3,4)}$ such that $\pi_4(\bGamma^{(3,4)})={E}$.

Let us summarize. The boundary between the second and third sheets goes along the family of~$p$
closed disjoint curves  $\bGamma^{(2,3)}$ such that $\pi_4(\bGamma^{(2,3)})=\pi_2(\FF)=F$, ${E}\cap F=\varnothing$.
Note that with this definition it is quite possible that both points
$\infty^{(2)},\infty^{(3)}$ lie in $\bGamma^{(2,3)}$.  Note also that the boundary
of the compact set~ $\FF$ has no common points with the compact set
$\bGamma=\partial\RS^{(1)}_2(w)=\partial\RS^{(2)}_2(w)$, and moreover, $\bGamma^{(2,3)}$
is disjoint both from  $\bGamma^{(1,2)}$ and from $\bGamma^{(3,4)}$ (cf.\ the assumption in~\cite{Kom20}).


It is clear that by performing the above procedure of continuation of the Green function
$g_{\FF}(\zz,\infty^{(1)})$ from the surface $\mathfrak
P_1\simeq\RS_2(w)\setminus\FF$ to the entire Riemann surface  $\RS_4(f)$ we get the bipolar
Green function $g(\zz,\infty^{(1)},\infty^{(4)})$ defined  on the Riemann surface
$\RS_4(f)$. This function has logarithmic singularities at the points
$\zz=\infty^{(1)}$ and $\zz=\infty^{(4)}$ and is normalized by the condition
$g(\zz,\infty^{(1)},\infty^{(4)})\equiv0$ for $\zz\in\bGamma^{(2,3)}$
(see \cite{Chi18},~\cite{Chi19}). The fact that the Nuttall partition  of a~Riemann surface
into sheets can be related  to the (zero, with appropriate normalization) level line
of the bipolar Green function for the Riemann surface  was known earlier only for the
case of Pad\'e polynomials  and the two-sheeted Riemann surface associated with the original germ.
In this classical case, the bipolar Green function spits the Riemann surface into two sheets.

It is easily checked that  the Green function $g(\zz):=g_{E}(z,\infty)$, $\zz=z^{(1)}\in\RS^{(1)}_2(w)$,
also extends to the four-sheeted Riemann surface $\RS_4(f)$.  Indeed,
$g(\zz)=g(\zz,\infty^{(1)},\infty^{(2)})$ is the bipolar Green function
for the Riemann surface $\RS_2(w)$. Since the points of the compact set~$\FF$ are not singular for this function,
one can define, for a~give partition of  $\RS_4(f)$ into sheets,
\begin{equation} \begin{aligned}
g(z^{(3)}):&=g(z^{(2)}),\qquad z\in\myh\CC,\\
g(z^{(4)});&=g(z^{(1)}),\qquad z\in\myh\CC\setminus{E} \end{aligned}
\label{11}
\end{equation}
with the given partition of the Riemann surface.
The function $g(\zz)$ thus obtained has logarithmic singularities
at all four\footnote{If $\infty\in F$, then $\infty^{(2)}=\infty^{(3)}$ and in this case, it is necessary, at this point at infinity
lying on the Riemann surface $\RS_4(f)$, to introduce the corresponding
local coordinate and then argue as in the paper \cite{KoPaSuCh17}.} points of the set $\pi_4^{-1}(\infty)$:
$\infty^{(1)}$, $\infty^{(2)}$, $\infty^{(3)}$,~$\infty^{(4)}$.

Now the following definition is correct.   For $\zz\in\RS_4(f)$, we put
\begin{equation}
u(\zz):=-2g_{\FF}(\zz,\infty^{(1)})-g(\zz),\quad
\zz\notin\pi^{-1}_4(\infty).  \label{12}
\end{equation}

The following result holds.

\begin{theorem}\label{the2} Let~$u(\zz)$, $\zz\in\RS_4(f)$, be the function defined by~\eqref{12}.

{\rm 1)} The following asymptotic formulas hold:
\begin{equation}
\begin{aligned} u(\zz)&=-3\log{|z|}+O(1),\quad\zz\to\infty^{(1)},\\
u(\zz)&=\log{|z|}+O(1) , \end{aligned} \label{13}
\end{equation}
 as $\zz$ tends to any of the points of the set
$\pi_4^{-1}(\infty)\setminus\infty^{(1)}=\{\infty^{(2)},\infty^{(3)},\infty^{(4)}\}$.

{\rm 2)} For the above partition of the Riemann surface $\RS_4(f)$ into sheets,
\begin{equation}
u(z^{(1)})<u(z^{(2)})<u(z^{(3)})<u(z^{(4)})\quad\text{for} \quad
z=\pi_4(z^{(j)})\notin E\cup F, \quad j=1,2,3,4.  \label{14}
\end{equation}
\end{theorem}

Together relations \eqref{13} and~\eqref{14} mean that the function
$u(\zz)$, as defined by~\eqref{12} and which is harmonic on
$\RS_4(f)\setminus\pi^{-1}_4(\infty)$, coincides up to a~constant
to the real part of the Abelian integral used in the definition
of the Nuttall partition of the Riemann surface  into sheets (see \cite{Nut84},~\cite{KoPaSuCh17}).

Thus, the partition of the Riemann surface $\RS_4(f)$ into sheets
$\RS_4^{(1)}(f)$, $\RS_4^{(2)}(f)$, $\RS_4^{(3)}(f)$, $\RS_4^{(4)}(f)$,
which was constructed with the use of the extremal problem~\ref{pro1}, is a~Nuttall partition.

\section{Proof of Theorem \ref{the1}}\label{s2}


So, given an $f\in\mathscr Z$ defined by \eqref{2}, suppose that the above assumptions
on the geometric arrangement of branching points and the corresponding  exponents are satisifed.

We consider the two-sheeted Riemann surface $\RS_2(w)$ of the function~$w$ defined by
$ w^2=z^2-1$.  We shall assume that the Riemann surface $\RS_2(w)$ is realized as a~two-sheeted
covering of the Riemann surface  $\myh\CC$ using the explicitly given uniformization:
\begin{equation} z=\frac12\left(\zeta+\frac1\zeta\right), \quad
w=\frac12\left(\zeta-\frac1\zeta\right),\quad \zeta\in\myh\CC_\zeta.  \label{32.0}
\end{equation}
Accordingly, the first sheet $\RS^{(1)}_2(w)$, on which $w=(z^2-1)^{1/2}/z\to1$ as $z\to\infty$, corresponds to the exterior
$\UU_\zeta:=\{\zeta:|\zeta|>1\}$ of the unit  disk
$\DD_\zeta:=\{\zeta:|\zeta|<1\}$ in the $\zeta$-plane, and the second sheet
$\RS^{(2)}_2(w)$, on which  $w=-(z^2-1)^{1/2}/z\to-1$ as $z\to\infty$,
corresponds to the unit disk~$\mathbb D_\zeta$ itself. By a~point~$\zz$ on the Riemann surface $\RS_2(w)$,
$\zz\in\RS_2(w)$, we shall mean the pair $\zz:=(z,w)=(z,\pm(z^2-1)^{1/2})$.
The canonical projection $\pi_2\colon\RS_2(w)\to\myh\CC$ is defined by $\pi_2(\zz):=z$.
By the point $\zz=\infty^{(1)}\in\RS_2(w)$
we mean
the point on the Riemann surface $\RS_2(w)$ such that $\pi_2(\infty^{(1)})=\infty$ and
$w/z\to1$ as $\zz\to\infty^{(1)}$. Similarly,, for $\zz=\infty^{(2)}$ we have  $\pi(\infty^{(2)})=\infty$ and $w/z\to-1$ as $\zz\to\infty^{(2)}$.
A~passage from the first sheet $\RS^{(1)}_2(w)$ to the second sheet $\RS^{(2)}_2(w)$
proceeds along the cut closed interval~${E}$. Here we assume as usual that the interval has to edges (the upper and the lower ones)
and that the sheets are ``glued'' by identifying crosswisely the edges
of the cuts; that is,
by identifying the upper edge of one cut with the lower edge of the
other cut, and vice versa.
It can be easily shown that the above partition into sheets is
a~Nuttall partition. Indeed, let
\begin{equation}
G(\zz):=\int_{-1}^{\zz}\frac{dt}{\sqrt{t^2-1}}=\log(z+w)=\log(z\pm(z^2-1)^{1/2})
\label{32}
\end{equation}
be an Abelian integral of the third kind with purely imaginary periods and logarithmic
singularities only at the points $\zz=\infty^{(1)}$ and $\zz=\infty^{(2)}$. Hence
$\eta_2(\zz):=-\Re G(\zz)$ is a~harmonic function on the Riemann surface
$\RS_2(w)\setminus\{\infty^{(1)}, \infty^{(2)}\}$,
$$
\eta_2(\zz)=\mp\log|z|+O(1),\quad \zz\to\infty^{(1)}\quad\text{or}\quad\zz\to\infty^{(2)},
$$
and
\begin{equation}
\eta_2(z^{(1)})<\eta_2(z^{(2)}),\quad z\in D.
\label{33}
\end{equation}
Moreover,  $\eta_2(\zz)=0$ for $\zz\in\bGamma$, $\bGamma:=\pi^{-1}_2({E})$.
So, $-\eta_2(z^{(1)})=g_{{E}}(z,\infty)$ is the Green function of the domain
$D$.  The given germ $f_\infty\in\HH(\infty)$ of a~function $f\in\SZ$, $f\in\HH(D)$, is lifted to the point
$\zz=\infty^{(1)}\in\RS_2(w)$ and extends to the entire first sheet
$\RS^{(1)}_2(w)$ as a~single-valued holomorphic  function (recall that
$\pi_2(\RS_2^{(1)}(w))=D=\myh\CC\setminus{E}$). A~further single-valued extension of this
function to the entire second sheet of the Riemann surface $\RS^{(2)}_2(w)$
is hindered by the branch points $a^{(2)}_j\in\RS^{(2)}_2(w)$ such that
$\pi_2(a^{(2)}_j)=a_j$, $j=1,\dots,2p$, $a_j=(A_j+1/A_j)/2$. In order that such a~single-valued
analytic  (meromorphic) extension of the germ $f_\infty$ be possible one
should make appropriate cuts on the second sheet $\RS_2^{(2)}(w)$.  Above in~\S~\ref{s1} we introduced
the corresponding class of admissible compact sets $\myk(f_\infty)$.


Since the Riemann surface $\RS_2(w)$ is of zero genus, Theorem~\ref{the1} can be reduced to the planar case
and to the corresponding compact set of minimal capacity on the plane. Indeed, the uniformization of $\RS_2(w)$
is defined using the Zhukovskii function~\eqref{32.0}.  Here, to the first sheet  of  $\RS^{(1)}_2(w)$
there corresponds the exterior $\UU:=\myh\CC\setminus\myo\DD$ of the unit disk
$\DD=\{\zeta:|\zeta|<1\}$, and to the second sheet, the unit disk~$\DD$. So, in this case, the quantity
$1/\zeta$ is the local coordinate~$\xi$ alluded to above. With this uniformization,
to each compact set $\KK\in\mathfrak
K(f_\infty)$ admissible for the germ~$f_\infty$  there corresponds an admissible compact set $\myt{\KK}\in\mathfrak
K(\myt{f}_\infty)$, $\myt{\KK}\subset\DD_\zeta$, for the germ
$\myt{f}_\infty\in\HH(\infty)$ corresponding to $f_\infty$. The Green function is invariant with respect to
conformal mappings of the domain. The Robin constant changes accordingly. Hence Problem~\ref{pro1} on the maximum
of the Robin constant $\gamma(\KK)$ on the class of compact sets $\KK\in\mathfrak
K(f_\infty)$ is equivalent to the maximization problem of the Robin constant $\gamma(\myt{K})$
on the class of admissible compact sets $\myt{K}\subset\DD_\zeta$,
$\myt{K}\in\mathfrak K(\myt{f}_\infty)$; that is, it is equivalent to the Stahl problem on
compact set of minimal capacity.

Under transformation~\eqref{32.0}, a~function $f(z)\in\SZ$ of the form~\eqref{2}
is transformed to the function
\begin{equation}
\myt{f}(\zeta)=\prod_{j=1}^{2p}\left(A_j-\frac1{z+w}\right)^{\alpha_j}
=\prod_{j=1}^{2p}\left(A_j-\frac1\zeta\right)^{\alpha_j},
\label{38}
\end{equation}
where all $\alpha_j=\pm1/2$. Since $A_j\in\UU$ for all~$j$, all singular points
of the function~$\myt{f}$ have the form  $\zeta=\myt{A}_j=1/A_j\in\DD$. Thus, to the extremal compact set~$\FF$
for the function~$f$ on the Riemann surface $\RS_2(w)$
there corresponds an admissible compact set of minimal capacity  $\myt{F}$ for the function
$\myt{f}$, which is given by $\myt{f}_\infty\in\HH(\infty)$ (note that
$\myt{f}\in\MM(\myh\CC\setminus\myt{F})$ and at the point $\zeta=0$ the function $\myt{f}$ has a~pole).
By the well-known properties of a~compact set of
minimal capacity, we have $\myt{F}\subset\DD$ (more precisely, the compact set $\myt{F}$
lies in the convex hull of the set
$\{\myt{A}_j=1/A_j,j=1,\dots,2p\}$) and $\myt{F}$ consists of a~finite number of
analytic arcs (which are trajectories of the quadratic differential), does not split the complex plane, and has
the  $S$-property~\eqref{7}.  Moreover (see \eqref{8}),
\begin{equation}
\myt{F}=\biggl\{\zeta\in\CC:\Re\int_{\myt{A}_1}^\zeta\sqrt{\frac{V_{2p-2}(t)} {B_{2p}(t)}}\,dt=0\biggr\},
\label{39}
\end{equation}
where $B_{2p}(t):=\prod_{j=1}^{2p}(t-\myt{A}_j)$, $V_{2p-2}(t):=(t-v_1)\dots(t-v_{2p-2})$
is the corresponding Chebotarev polynomial, $v_j$, $j=1,\dots,2p-2$, are the Chebotarev   points
of the compact set~$\myt{F}$. All these properties  of  a~compact set of minimal capacity are
well known in the general case and have been obtained already by Stahl in~1985
(see \cite{Sta97b}, ~\cite{Sta12}, and also~\cite{ApBuMaSu11}).

In the case considered here all the branch points of the function~$\myt{f}$
are of second order. In this setting, the existence and description of a~compact set of
minimal capacity was given already by Nuttall~\cite{NuSi77},~\cite{Nut84}. In particular,
in the case of general position, all the zeros of the polynomial $V_{2p-2}$ are of even multiplicity, all $v_j\neq
\myt{A}_k$ for  $j=1,\dots,2p-2$, $k=1,\dots,2p$, and the compact set $\myt{F}$ consists of
$p$~disjoint analytic arcs that pairwisely connect points of the set
$\{\myt{A_j},j=1,\dots,2p\}$. Moreover, the compact set~$\myt{F}$ has the classical $S$-property; that is, the
the corresponding
Green function satisfies the relation of the form~\eqref{7}.
By the invariance of the Green function with respect to conformal mappings
it follows that the Green function corresponding to the compact set~$\FF$ satisfies~\eqref{6}. Hence the
compact set $\FF$ has the required $S$-property.

So, on the second sheet of the Riemann surface $\RS_2(w)$ we have obtained a~system of analytic arcs comprising the compact set $\FF$,
not splitting the Riemann surface $\RS_2(w)$, having the
$S$-property~\eqref{6}, and such that
$f_\infty\in\MM(\RS_2(w)\setminus{\FF})$. The $S$-property~\eqref{6} is a~direct consequence
of the $S$-property from~\eqref{7}, because the Green function
is invariant with respect to a~conformal mapping.

\section{Proof of Theorem~\ref{the2}}\label{s3}

The first part of Theorem~\ref{the2} is a direct consequence of
Definition~\eqref{12} of the function~$u(\zz)$ and properties of the functions $g(\zz)$ and
$g_{\FF}(\zz,\infty^{(1)})$.

Let us now prove the second part of Theorem~\ref{the2}.


We proceed as follows. Consider the function
\begin{equation*}
u(\zz):=-2g_{\FF}(\zz,\infty^{(1)})-g(\zz),\quad \zz\in\RS_2(w)\setminus\FF,
\end{equation*} where $g(\zz)=\log|z+w|=\log|\zeta|$. The function $u(\zz)$ is harmonic  in the domain $\RS_2(w)\setminus\FF$,
with the exception of
the points at infinity  $\zz=\infty^{(1)}$ and $\zz=\infty^{(2)}$, where it behaves as follows:
\begin{equation}
u(\zz)= \begin{cases}
-3\log|z|+O(1),&\zz\to\infty^{(1)},\\ \log|z|+O(1),&\zz\to\infty^{(2)}
\end{cases} \label{40}
\end{equation}
 (here and in what follows we assume for simplicity that $\infty^{(2)}\not\in\FF$;
 otherwise a~more careful consideration is required quite similar to that conducted in~\cite{KoPaSuCh17}).

We set, as before
 \begin{equation}
\bGamma=\bigl\{\zz\in\RS_2(w):u(z^{(1)})=u(z^{(2)})\bigr\};
\label{41}
\end{equation}
here when writing $z^{(1)}\in\RS^{(1)}_2(w)$ and
$z^{(2)}\in\RS^{(2)}_2(w)$ we mean, as in the above, a~Nuttall partition of the Riemann surface $\RS_2(w)$ into sheets, $\pi_2(\bGamma)=E$.
The compact set~$\bGamma$ is a~closed arc on the Riemann surface $\RS_2(w)$ passing through the points $\zz=\pm1$,
not intersecting the compact set $\FF$, and splitting $\RS_2(w)$ into two domains, of which one contains~$\FF$,
and the other one, the point  $\zz=\infty^{(1)}$.


For  $z\notin E\cup F$, we set
 $$
 v_1(z):=u(z^{(2)})-u(z^{(1)}).
 $$
 By \eqref{12}, using properties of the functions $g(\zz)$ and $g_{\FF}(\zz,\infty^{(1)})$ and taking into account the
symmetry of the compact set~$F$ with respect to the real line, we have
\begin{align*}
v_1(z)&=-2g_{\mathbf F^{(2)}}(z^{(2)},\infty^{(1)})-g(z^{(2)})+2g_{\mathbf F^{(2)}}(z^{(1)},\infty^{(1)})+g(z^{(1)}) \\
&=2g_{\mathbf F^{(2)}}(z^{(1)},\infty^{(1)})-2g_{\mathbf F^{(2)}}(z^{(2)},\infty^{(1)}) +2g(z^{(1)}).
\end{align*}
This implies the following properties  of the function~$v_1$:

1) $v_1(z)$ is a~harmonic function in the domain $\CC\setminus(E\cup F)$ and is continuous in~$\CC$;

2) $v_1(z)\equiv0$ for $z\in E$;

3) $v_1(z)=4\log|z|+O(1)$ as $z\to\infty$, $z\notin F$;

4) $v_1(z)>0$ for $z\in F\setminus\infty$.

From properties 1)--4) of the function~$v_1$ we get the inequality $v_1(z)>0$ for
$z\in\CC\setminus E$, and therefore, the inequality
$$
u(z^{(1)})<u(z^{(2)})\quad\text{for}\quad z\in\myh{\CC}\setminus E.
$$
We set $v_2(z):=u(z^{(3)})-u(z^{(2)})$, $z\notin E\cup F$.
By definition \eqref{12} and using the properties of the function~$g(\zz)$
we have $v_2(z)=-2g_{\FF}(z^{(3)},\infty^{(1)})+2g_{\FF}(z^{(2)},\infty^{(1)})=4g_{\FF}(z^{(2)},\infty^{(1)})>0$.
Therefore,
$$
u(z^{(2)})<u(z^{(3)})\quad\text{for}\quad z\in\myh{\CC}\setminus(E\cup
F).
$$

Finally, let us prove that
\begin{equation}
u(z^{(3)})<u(z^{(4)})\quad\text{for}\quad z\in\myh{\CC}\setminus(E\cup F).
\label{15}
\end{equation}
Indeed, we set
\begin{align*}
v_3(z):&=u(z^{(4)})-u(z^{(3)}) \\
&=2g_{\FF}(z^{(3)},\infty^{(1)})+g(z^{(3)})-2g_{\FF}(z^{(4)},\infty^{(1)})-g(z^{(4)}).
\end{align*}
By definition we have
\begin{gather*}
g_{\FF}(z^{(3)},\infty^{(1)})=-g_{\FF}(z^{(2)},\infty^{(1)})\quad\text{for}\quad z\notin F, \\
g_{\FF}(z^{(4)} ,\infty^{(1)})=-g_{\FF}(z^{(1)},\infty^{(1)})\quad\text{for}\quad z\notin F, \\
g(z^{(3)})=g(z^{(2)}),\qquad
g(z^{(4)})=g(z^{(1)}).
\end{gather*}
Hence
\begin{align*}
v_3(z)&=u(z^{(4)})-u(z^{(3)}) \\
&=-2g_{\FF}(z^{(2)},\infty^{(1)})
+2g_{\FF}(z^{(1)},\infty^{(1)})+g(z^{(2)})-g(z^{(1)}) \\
&=2g_{\FF}(z^{(1)},\infty^{(1)})-2g_{\FF}(z^{(2)},\infty^{(1)})
-2g(z^{(1)}).
\end{align*}
Hence the function $v_3(z)$, $z\in\myh{\CC}\setminus(E\cup F)$, has the following
properties:

1) the function $v_3$ is harmonic in the domain $\CC\setminus(E\cup F)$ and
is continuous  on the plane $\CC$;

2) $v_3(z)\equiv0$ for $z\in E$;

3) $v_3(z)=2g_{\FF}(z^{(1)},\infty^{(1)})-2g(z^{(1)})>0$ for $z\in F$;

4) for $\infty\notin F$, the function $v_3$ is harmonic near the point
$z=\infty$; if $\infty\in F$, then $v_3$ is a~continuous function at the point
$\infty\in F$ and is harmonic in $U\setminus F$, where $U$~is some
neighborhood of the  point $z=\infty$.

Indeed, properties  1), 2) and 4) of the function~$v_3$ are clear. Let us verify
property 3). Consider the  function
$$
v_4(z):=2g_{\FF}(z^{(1)},\infty^{(1)})-2g(z^{(1)}),\quad z\in\CC\setminus
E.
$$
We have $v_4(z)=v_3(z)$ for $z\in F$. The function $v_4(z)$ is continuous  on~$\CC$ and by
properties of the functions $g_{\FF}(\zz,\infty^{(1)})$ and $g(\zz)$
it extends as a~harmonic to the neighborhood of the  point $z=\infty$. The compact set~$\bGamma$~is an
admissible set\footnote{More precisely, this property is satisfied for the compact set
$\myo{\RS}_2^{(2)}(w)\supset\bGamma$, which is the closure of the second sheet of the Riemann surface $\RS_2(w)$ of the function
$w^2=z^2-1$.} for Problem~\ref{pro1}.  Since the compact set
$\FF$ lies on the second sheet of the Riemann surface  $\RS_2(w)$, $\FF\subset\RS_2^{(2)}(w)$,
we have  $g_{\FF}(\zz,\infty^{(1)})>0$ for $\zz\in\EE$ and
$g(\zz)\equiv0$ for $\zz\in\bGamma$.
It is clear that  for the function~$v_4(z)$ we have $v_4(z)>0$ for $z\in E$ and
$v_4(z)$~is a~harmonic function in $\myh{\CC}\setminus E$. Therefore,
$v_4(z)>0$ everywhere in $\myh{\CC}$, and hence, also on the compact set~$F$. So,
for $z\in F$, we have $v_3(z)=v_4(z)>0$. This proves property~3) of the function~$v_3$. Hence
$v_3(z)>0$ for $z\notin E\cup F$.
Now the required inequality
$u(z^{(3)})<u(z^{(4)})$ for $z\notin E\cap F$
follows from the definition of the function~$v_3(z)$.

So, for the function~$u(\zz)=-2g_{\FF}(\zz,\infty^{(1)})-g(\zz)$ we see that  $u(\zz)$ is a~harmonic function in the domain
$\RS_4(f)\setminus\{\infty^{(1)},\infty^{(2)},\infty^{(3)},\infty^{(4)}\}$.
Moreover, by the above we have
\begin{equation*}
u(\zz)= \begin{cases} -3\log|z|+O(1),&\zz\to\infty^{(1)},\\
\log|z|+O(1),&\zz\to\infty^{(2)},\\ \log|z|+O(1),&\zz\to\infty^{(3)},\\
\log|z|+O(1),&\zz\to\infty^{(4)},
\end{cases}
\end{equation*}
and further, under the above partition of the Riemann surface $\RS_4(f)$ into sheets, the following inequalities hold:
\begin{equation}
u(z^{(1)})<u(z^{(2)})<u(z^{(3)})<u(z^{(4)}).  \label{43}
\end{equation}

So, the partition of the Riemann surface  $\RS_4(f)$ into sheets, which we have introduced with the help of the Green
function  $g_{\FF}(\zz,\infty^{(1)})$ corresponding to the  extremal compact set of minimal capacity
$\FF\subset\RS_2(w)$, is shown to be a~Nuttall partition.

Theorem~\ref{the2} is proved.

\section{Some concluding remarks}\label{s4}

\subsection{}\label{s4s1}

\begin{figure}
\includegraphics[width=\textwidth]{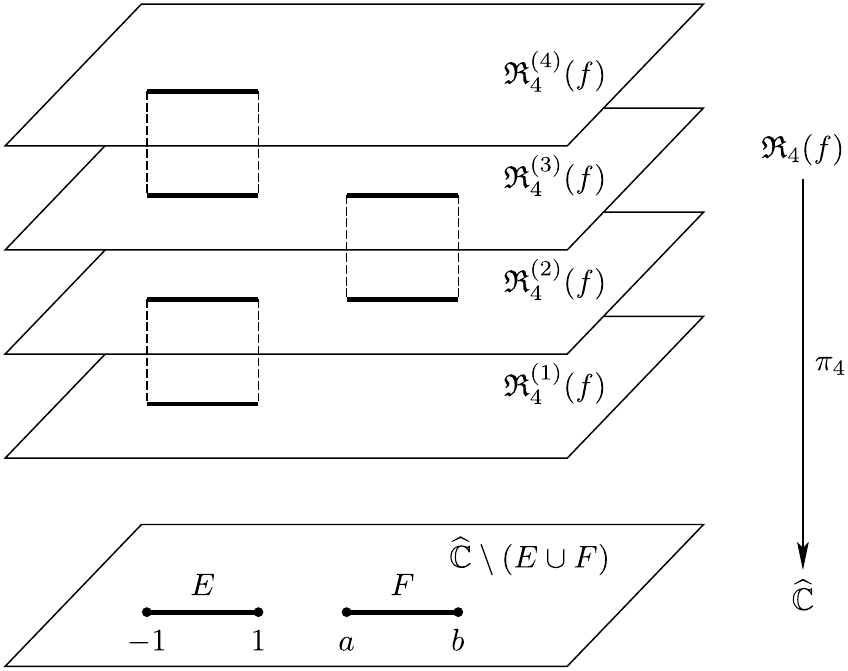}
\caption{Four-sheeted Riemann surface  $\RS_4(f)$ of the function~$f$ of the form~\eqref{1}.}
\label{fig_1}
\end{figure}

\begin{figure}
\includegraphics[width=\textwidth]{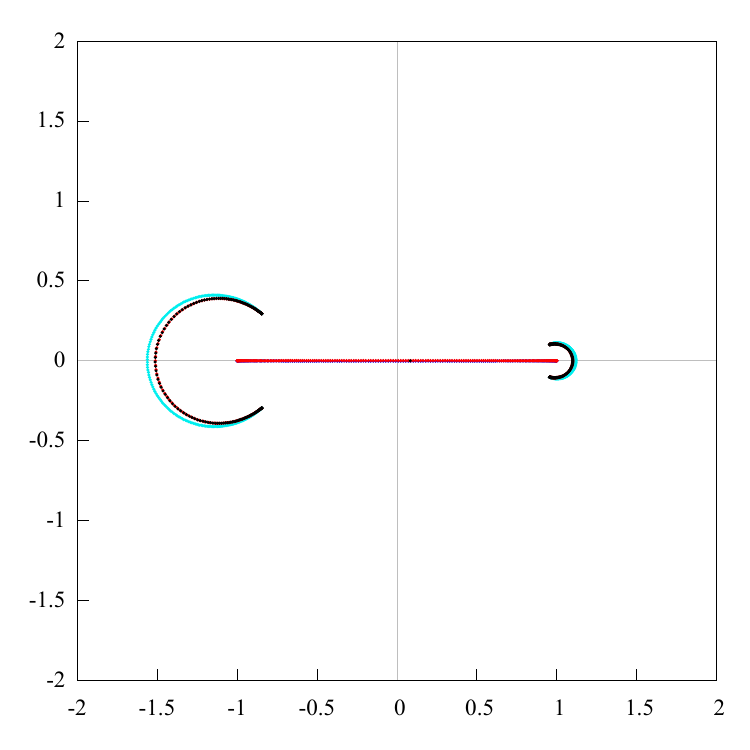}
\caption{}
\label{fig_2}
\end{figure}

In is natural, while remaining in the same class of multivalued functions~$\mathscr Z$,
to compare the solution of Problem~\ref{pro1}, as obtained in terms of the Nuttall partition
of the four-sheeted  Riemann surface  $\RS_4(f)$ of a~function $f\in\mathscr Z$, and
the existence of a~three-sheeted Riemann surface with a~Nuttall partition $\RS_3(f_\infty)$
associated with a~given germ  $f_\infty\in\HH(\infty)$; see \cite{Nut84},~\cite{RaSu13},~\cite{Sue19}.
The latter problem  will be referred to as Problem~\ref{pro2}; the solutions of these two
problems will be denoted by $\bGamma_4^{(2,3)}$
(above in the present paper it was denoted by $\mathbf F$) and
$\bGamma_3^{(2,3)}$, respectively.
The natural question is whether these compact sets can be equal.
More precisely, the question of course should be put like this: may the canonical projections
$\Gamma_4^{(2,3)}:=\pi_4(\bGamma_4^{(2,3)})$ and
$\Gamma_3^{(2,3)}:=\pi_3(\bGamma_3^{(2,3)})$  coincide with each other?
It is easily seen that this is indeed so in the  real-case situation, that is, when in \eqref{2} all (pairwise different)
quantities $A_j$ lie in $\RR$. In the case when
in $p=1$ and $1<A_1<A_2$ in \eqref{2}, the pair of functions $f,f^2$ forms a~Nikishin system  (see \cite{Sue18}).
The  four-sheeted Riemann surface  corresponding to this function~$f$ is depicted in Fig.~\ref{fig_1}. In this case,
$\pi_4(\bGamma_4^{(2,3)}) =\pi_3(\bGamma_3^{(2,3)})=[a_1,a_2]$, where
$a_j=(A_j+1/A_j)/2$. The question
of whether
$\pi_4(\bGamma_4^{(2,3)})$ and $\pi_3(\bGamma_3^{(2,3)})$ may coincide for a~function
$f\in\mathscr Z$ in some other (nonreal) case remains open from the theoretical point of view.
Numerical experiments surely show that
such a coincidence should not be expected if the reality condition is violated  (see Fig.~\ref{fig_2}).

We recall (see \cite{Sue19},~\cite{IkSu20}) that the  compact set $\Gamma_3^{(2,3)}$
is the solution of the following potential theory
``$\max$-$\min$''-problem for the nonstandard potential with harmonic external field.

\begin{problem}\label{pro2}
Let $\mathfrak K(f_\infty)$ be the above family of admissible compact sets  $\KK=\mathbf K^{(2)}\subset\RS_2^{(2)}$ for the germ
$f_\infty$, that is,

1) the set $\Omega(\KK):=\RS_2(w)\setminus\KK\ni\infty^{(1)}$ is a~domain  on the Riemann surface $\RS_2(w)$ of the function $w^2=z^2-1$;

2) $f\in\MM(\Omega(\KK))$.

Let $\bmu\in M_1(\KK)$ be a~unit Borel measure supported in the compact set $\KK$. We set  $\Phi(\zz):=z+w$.
Let
\begin{equation}
P^{\bmu}(\zz):=\int_{\KK}\log\frac{|1-1/(\Phi(\zz)\Phi(\mytt))|}{|z-t|^2}\,d\bmu(\mytt),
\label{4.1}
\end{equation}
be the potential of the measure $\bmu$, and let
$$
J_V(\bmu):=\int_{\KK}P^{\bmu}(\zz)\,d\dmu(\zz)+2\int_{\KK}V(\zz)\,d\dmu(\zz)
$$
be the energy of the measure  $\bmu$ corresponding to this potential and the external field $V(\zz):=-\log|\Phi(\zz)|$.

Let a measure $\blambda_{\KK}\in M_1(\KK)$ be the solution of the extremal problem
$$
J_V(\blambda_{\KK})=\min_{\bmu\in
M_1(\KK)}J_V(\bmu).
$$
Then  (see \cite{Sue19}) the ``$\max$-$\min$-problem''
\begin{equation}
J_V(\blambda_{\myt\FF})=\max_{\KK\in\mathfrak K(f_\infty)}J_V(\blambda_{\KK})
=\max_{\KK\in\mathfrak K(f_\infty)}\min_{\bmu\in M_1(\KK)}J_V(\bmu)
\label{4.3}
\end{equation}
has a~unique (in the class $\mathfrak K(f_\infty)$) solution
$\myt{\FF}\in\mathfrak K(f)$ in the class of compact sets $\KK\in\mathfrak K(f_\infty)$. The extremal compact set $\myt{\FF}$ has the following  $S$-property:
\begin{equation}
 \frac{\partial
\bigl(P^{\blambda_{\myt\FF}}(\zz)+V(\zz)\bigr)}{\partial n^{+}}
=\frac{\partial \bigl(P^{\blambda_{\myt\FF}}(\zz)+V(\zz)\bigr)}{\partial
n^{-}},\quad \zz\in\myt\FF^\circ;
\label{4.4}
\end{equation}
here $\myt\FF^\circ$  is the family of open arcs whose closures constitute~$\FF$.
\end{problem}

Moreover  (see \eqref{4.10}),
$$
\frac1n\chi^{\vphantom{p}}_{Q_{n,j}}\overset{*}\longrightarrow\pi_2(\blambda_{\FF}),
\quad n\to\infty.
$$
where the measure  $\pi_2(\blambda_{\FF})$ is defined   as $\pi_2(\blambda_{\FF})(e)=
\blambda_{\FF}(\mathbf e)$ for any $\mathbf e\subset\RS^{(2)}_2(w)$, $e=\pi_2(\mathbf e)$.

Note that the extremal measure $\blambda_{\KK}$ is a~(unique) equilibrium measure; that is,
\begin{equation}
P^{\blambda_{\KK}}(\zz)+V(\zz)\equiv w_{\KK}=\const,\quad \zz\in\KK.
\label{4.7}
\end{equation}
The facts that in \eqref{4.7} the identity
holds on the entire compact set $\KK$ and $\supp\blambda_{\KK}=\KK$ were proved in~\cite{Sue19b}.

The extremal problems~\ref{pro1} and~\eqref{4.3} are different.
It is natural to assume that the corresponding extreme
compact sets are also different.
More precisely, one may assume that in general
$\pi_3(\bGamma^{(2,3)}_3)\neq\pi_4(\bGamma^{(2,3)}_4)$, where
$\bGamma_3^{(2,3)}=\myt\FF$, $\bGamma_4^{(2,3)}=\FF$.
The fact that in some cases these two sets coincide follows from the example presented in Fig.~\ref{fig_1}.
This case corresponds to the choice of the
parameters $p=1$ and $1<A_1<A_2$ in~\eqref{2}. With this choice,
$\pi_3(\bGamma^{(2,3)}_3)=\pi_4(\bGamma^{(2,3)}_4)=[a_1,a_2]$.
At present, in the general case
we can only resort to numerical experiments.

For example, in Fig.~\ref{fig_2} we show the zeros (dark blue, red and black points) of three
Hermite--Pad\'e polynomials of type~I\enskip $Q_{300,j}$, $j=0,1,2$ (for the family $[1,f,f^2]$, see \eqref{4.10}).
From the calculated eight Hermite--Pad\'e polynomials of type~I\enskip
$q_{300,j}$ and $q_{299,j}$, $j=0,1,2,3$ (for the  family $[1,f,f^2,f^3]$, see \eqref{4.11}), we evaluated new
(nonstandard) Hermite--Pad\'e polynomials, which were introduced in~\cite{Sue18d} (see also~\cite{Kom20}).
The zeros of these new polynomials localize the projection
$\pi_4(\bGamma^{(2,3)}_4)$ on the Riemann sphere  $\myh{\CC}$ of the compact set
$\bGamma^{(2,3)}_4$, which is the boundary between the second and third Nuttall sheets
of the four-sheeted Riemann surface $\RS_4(f)$  of a~function $f\in\mathscr Z$.  The compact set
$\bGamma^{(2,3)}_4=\FF$ is the solution of problem~\ref{pro1}.  In Fig.~\ref{fig_2}, these zeros
of the new  (nonstandard) Hermite--Pad\'e polynomials
are shown by pale blue points. It is clear that the compact sets
$\pi_3(\bGamma_3^{(2,3)})$ and $\pi_4(\bGamma^{(2,3)}_4)$ differ from each other.
The red points located on the real line are the zeros
of the Pad\'e polynomials  of order $100$. They correspond to the closed interval $[-1,1]=:\Delta$.
Correspondingly, this is the projection  $\pi_2(\bGamma^{(1,2)}_2)$ of the boundary
$\bGamma^{(1,2)}_2$ between the first and second sheets of the Riemann surface $\RS_2(w)$ of the function
$w^2=z^2-1$.

\subsection{}\label{s4s2} Let us now consider the case of two intervals. Namely, let
 $\Delta_1=[e_1,e_2]$ and $\Delta_2=[e_3,e_4]$, where $e_1<e_2<e_3<e_4$, and let
$\pfi_{\Delta_j}$, $j=1,2$, be the functions inverse to the Zhukovskii function  corresponding to these intervals,
$|\pfi_{\Delta_j}(z)|>1$ for $z\notin\Delta_j$. We set
\begin{equation}
f(z)=\prod_{j=1}^{2p}\left(A_j-\frac1{\pfi_{\Delta_1}(z)}\right)^{\alpha_j}\prod_{k=1}^{2q}\left(B_k-\frac1{\pfi_{\Delta_2}(z)}\right)^{\beta_k},
\label{4.8}
\end{equation}
where   $|A_j|>1$, $|B_k|>1$, $A_j=\myo{A}_s$ for all~$j$ and for some $s\in\{1,\dots,2p\}$,
$B_k=\myo{B}_\ell$ for some $\ell\in\{1,\dots,2q\}$ and  all $A_k,B_j\notin\RR$,
$\alpha_j,\beta_k\in\{1/2,-1/2\}$,
$\sum_{j=1}^{2p}\alpha_j+\sum_{k=1}^{2q}\beta_k=0\pmod{\ZZ}$.  We denote by $\mathscr Z(\Delta_1,\Delta_2)$
the class of functions of the form~\eqref{4.8} with the above conditions on
$A_k$ and~$B_j$. By the assumption,
$|\pfi_{\Delta_j}(z)|>1$ for $z\notin\Delta_j$. Hence $f\in\HH(D)$, where
$D:=\myh{\CC}\setminus(\Delta_1\cup\Delta_2)$, and in particular,
$f\in\HH(\infty)$.  The set of branching points of the  multivalued function~$f$
consists of the  points $\pm1$ and the points $a_j=\pfi^{-1}_{\Delta_1}(A_j)\notin\RR$, $j=1,\dots,2p$,
$b_k=\pfi^{-1}_{\Delta_2}(B_k)\notin\RR$, $k=1,\dots,2q$.  It is clear that,  for $f\in\mathscr
Z(\Delta_1,\Delta_2)$, the Stahl compact set  consists of two closed intervals,
$S=\Delta_1\cup\Delta_2$, and  $D$~is the corresponding Stahl domain.

A natural question arises: in what terms should the limit
distribution of the zeros of the  Hermite--Pad\'e polynomials of type~I\enskip $Q_{n,j}$
(see \eqref{4.10}) for the family $[1,f,f^2]$ of functions~$f$ from the class
$\mathscr Z(\Delta_1,\Delta_2)$ be characterized?

The conjecture is that instead  of the Riemann surface of the function $w^2=z^2-1$ we should now
consider the Riemann surface of the function $w^2=(z-e_1)(z-e_2)(z-e_3)(z-e_4)$.
In this case, the general
form of the potential \eqref{4.1}, the external field, and the energy will be preserved.
The $S$-compact set~$\FF$, which corresponds to the problem of the limit distribution of the zeros of the  Hermite--Pad\'e polynomials
for the function  $f\in\mathscr Z(\Delta_1,\Delta_2)$,
is as before is characterized as the solution of an extremal problem
of the form \eqref{4.3}, while the $S$-property itself has the form~\eqref{4.4}.

It is worth pointing out that in \eqref{4.1}--\eqref{4.3} we speak about the $S$-compact set
$\myt{\mathbf F}^{(2,3)}$, whose existence is related to the existence
of a~three-sheeted Nuttall-par\-ti\-tion\-ed Riemann surface  $\RS_3(f_\infty)$ associated with a~given germ
$f_\infty\in\HH(\infty)$, $f\in\mathscr Z(\Delta_1,\Delta_2)$.
It is an open question whether there exists a~four-sheeted Nuttall-par\-ti\-tion\-ed Riemann surface
 $\RS_4(f_\infty)$ associated with a~germ
$f_\infty\in\HH(\infty)$. In the case $f\in\mathscr
Z(\Delta_1,\Delta_2)$ the Riemann surface  of the function~$f$ is 8-sheeted and the existing
solution of Problem~\ref{pro1}
apparently has nothing to do with the compact set
$\bGamma_4^{(2,3)}$, which is the boundary between the second and third sheets of the Riemann surface
$\RS_4(f_\infty)$. Note that since the parameters
$A_j$ and $B_k$ are real symmetric, the projection of the boundary $\bGamma_3^{(1,2)}$ between the first and second sheets
of the Riemann surface $\RS_4(f_\infty)$  always coincides with the union of the closed intervals~$\Delta_1$ and
$\Delta_2$. The same result also holds for the projection of the boundary $\bGamma_4^{(1,2)}$ between the first and second sheets of the Riemann surface~$\RS_4(f)$.

\begin{figure}
\includegraphics[width=\textwidth]{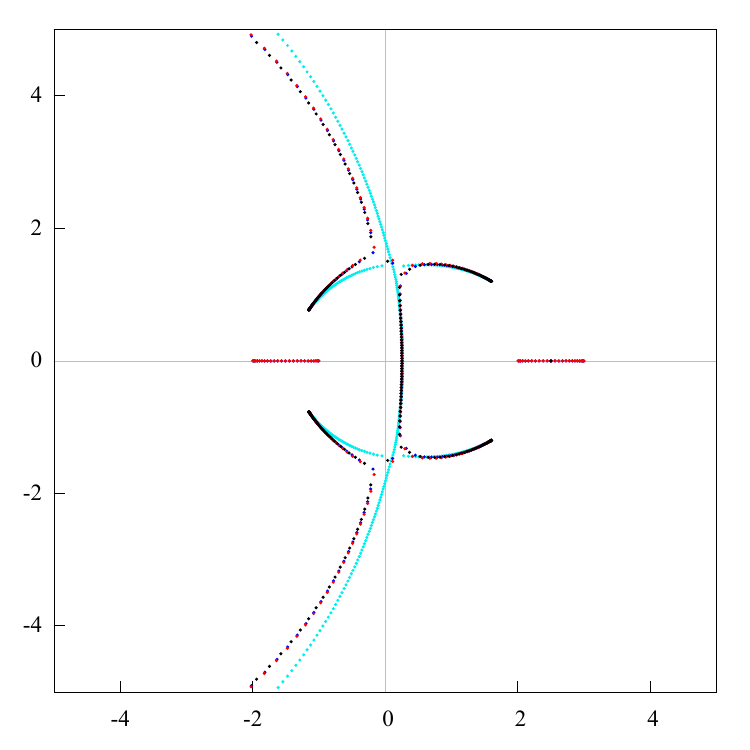}
\caption{}
\label{fig_3}
\end{figure}

Now   let  $p=q=1$, $\alpha_1=\alpha2=1/2$, $\beta_1=\beta_2=-1/2$ in \eqref{4.8}, that is,
\begin{equation}
f(z)=\left(\prod_{j=1}^2\left(A_j-\frac1{\pfi_{\Delta_1}(z)}\right)/\prod_{k=1}^2\left(B_k-\frac1{\pfi_{\Delta_2}(z)}\right)\right)^{1/2}.
\label{4.9}
\end{equation}
For a~function~$f$ of the form~\eqref{4.9}, under certain values of  $A_j$
and $B_k$ satisfying the above conditions,  we have numerically found the zeros
of the Hermite--Pad\'e polynomials  of type~I\enskip  $Q_{n,0},Q_{n,1},Q_{n,2}$, $\mdeg{Q_{n,j}}=n$,
$n=300$, for the family $[1,f,f^2]$ satisfying the relation
\begin{equation}
\left(Q_{n,0}+Q_{n,1}f+Q_{n,2}f^2\right)(z)=O\left(z^{-2n-2}\right), \quad
z\to\infty,
\label{4.10}
\end{equation}
and identified the zeros of the Hermite--Pad\'e polynomials
of type~I\enskip $q_{n,0}$, $q_{n,1}$, $q_{n,2}$, $q_{n,3}$, $\mdeg{q_{n,j}}=n$, $n=300$,  for the family
$[1,f,f^2,f^3]$ satisfying the relation
\begin{equation}
\left(q_{n,0}+q_{n,1}f+q_{n,2}f^2+q_{n,3}f^3\right)(z)=O\left(z^{-3n-3}\right),\quad
z\to\infty.
\label{4.11}
\end{equation}
According to the available theoretical results
and conjectures  (see \cite{Nut84},~\cite{RaSu13},~\cite{Sue19}), the zeros of the polynomials $q_{300,j}$
localize the compact set $\pi_3(\bGamma_3^{(2,3)})$, which is the projection
of the compact set  $\bGamma_3^{(2,3)} $ lying on the Riemann surface $\RS_3(f_\infty)$ onto the Riemann sphere $\myh{\CC}$.
Here,
$\RS_3(f_\infty)$ is a~three-sheeted Riemann surface with Nuttall partition into sheets associated with the germ
$f_\infty$, $\bGamma_3^{(2,3)}$ is the boundary between its
second and third sheets. In Fig.~\ref{fig_3}, these zeros of the polynomials $Q_{300,j}$,
$j=0,1,2$, are shown by dark blue, red, and black points.
It is clearly seen that the compact set  $\pi_3(\bGamma^{(2,3)})$ has four Chebotarev  points of zero density.

From the calculated polynomials of type~I\enskip $q_{300,j}$ and $q_{299,j}$, $j=0,1,2,3$,
we calculate new (nonstandard) Hermite--Pad\'e polynomials, which were introduced
in~\cite{Sue18d} (see also~\cite{Kom20}). The zeros of these polynomials localize the
projection $\pi_4(\bGamma^{(2,3)}_4)$ on the Riemann sphere $\myh{\CC}$ of the compact set
$\bGamma^{(2,3)}_4$, which  is the boundary between the second and third Nuttall sheets of
the four-sheeted Riemann surface $\RS_4(f_\infty)$ associated with the germ
$f_\infty\in\HH(\infty)$ of a~function $f\in\mathscr Z(\Delta_1,\Delta_2)$.  In Fig.~\ref{fig_3},
the zeros of the new (nonstandard) Hermite--Pad\'e polynomials are shown by pale blue points.
It is clearly seen that, first,  the compact set $\pi_4(\bGamma^{(2,3)}_4)$ contains four Chebotarev  points,
but they all have positive density. Second, the compact sets
$\pi_3(\bGamma_3^{(2,3)})$ and $\pi_4(\bGamma^{(2,3)}_4)$ differ from each other.

Next, the red points on the real line are the zeros
of the Pad\'e polynomials  of order $100$. They correspond to two closed intervals $\Delta_1$ and
$\Delta_2$. 
The set  $\Delta_1\cup\Delta_2$ - is the projection   $\pi_2(\bGamma^{(1,2)}_2)$of the boundary
between the first and second sheets of the Riemann surface   $\RS_2(w)$  of the function 
$w^2=(z-e_1)(z-e_2)(z-e_3)(z-e_4)$.

\begin{figure}
\includegraphics[width=\textwidth]{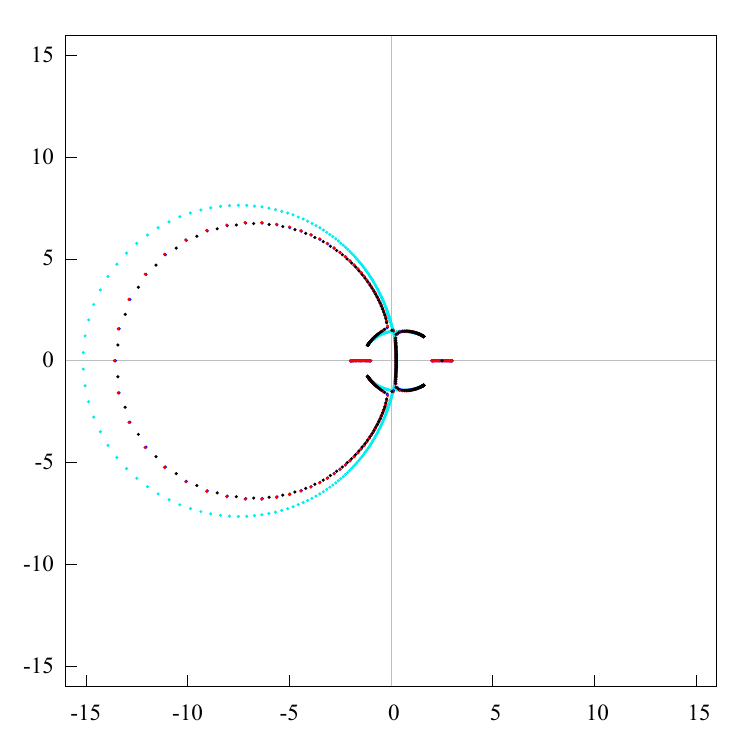}
\caption{}
\label{fig_4}
\end{figure}

In Fig.~\ref{fig_4} we show the same three sets, as in Fig.~\ref{fig_3}, but on a~smaller scale/
Each of the sets
$\pi_3(\bGamma^{(2,3)}_3)$ and $\pi_4(\bGamma^{(2,3)}_4)$ splits
the Riemann sphere  $\myh{\CC} $ into two domains so that each of these domains
contains precisely one closed interval $\Delta_1$ or $\Delta_2$. At the same time,
if  on considers the sets $\bGamma^{(2,3)}_3$ and
$\bGamma^{(2,3)}_4$ on the Riemann surface $\RS_2(w)$ of the function
$w^2=(z-e_1)(z-e_2)(z-e_3)(z-e_4)$, then it appears that the complement of each of such sets
is not a~domain on this Riemann surface.
This empirical fact should prove to be of utmost importance in generalizing the results obtained in~\cite{Lop20}
for a~(real) Nikishin system  (that is, for the case when the compact set~$F$ is the union of real closed intervals)
to the more general complex case, in which the compact set~$F$ is not known \textit{a~priori}; cf.~\cite{RaSu13},
where the case of a~complex Nikishin system is considered, but in the case when the complement of
$\Gamma_3^{(2,3)}= \pi_3(\bGamma_3^{(2,3)})$ is connected.

\clearpage
\newpage

\clearpage
\newpage

\clearpage
\newpage

\clearpage
\newpage

\clearpage
\newpage


\begin{thebibliography}{99}

\bibitem{Apt08} 
\href{http://www.ams.org/mathscinet-getitem?mr=2475084}
{A. I. Aptekarev, ``Asymptotics of Hermite--Pad\'{e} approximants for a pair of functions with branch points'',
\textit{Dokl. Math.}, \textbf{78}:2 (2008), 717--719.}

\bibitem{ApBuMaSu11} 
\href{http://mi.mathnet.ru/eng/umn9448}
{A. I. Aptekarev, V. I. Buslaev, A. Mart{\'\i}nez-Finkelshtein,
S. P. Suetin, ``Pad\'{e} approximants, continued fractions, and orthogonal polynomials'',
\textit{Uspekhi Mat. Nauk}, \textbf{66}:6(402) (2011), 37--122.}

\bibitem{ApYa15} 
\href{http://dx.doi.org/10.1007/s11511-016-0133-5}
{Alexander I. Aptekarev, Maxim L. Yattselev, ``Pad\'{e} approximants for functions with
branch points -- strong asymptotics of Nuttall--Stahl polynomials'', \textit{Acta Mathematica},
\textbf{215}:2 (2015), 217--280.}

\bibitem{ApTu16} 
\href{http://mi.mathnet.ru/eng/izv8420}
{A. I. Aptekarev, D. N. Tulyakov, ``Nuttall's Abelian integral
on the Riemann surface of the cube root of a polynomial of degree 3'', \textit{Izv. Math.},
\textbf{80}:6 (2016), 997--1034.}

\bibitem{Ara03} 
\href{http://www.ams.org/mathscinet-getitem?mr=2133957}
{N. U. Arakelyan, ``Efficient analytic continuation of power series and the localization
of their singularities'', \textit{Izv. Nats. Akad. Nauk Armenii Mat.}, \textbf{38}:4 (2003), 5--24;
in J. Contemp. Math. Anal. \textbf{38}:4 (2004), 2--20.}

\bibitem{Chi18} 
\href{http://mi.mathnet.ru/eng/tm3918}
{E. M. Chirka, ``Potentials on a compact Riemann surface'', \textit{Proc. Steklov Inst. Math.}, \textbf{301} (2018), 272--303.}

\bibitem{Chi19} 
\href{http://mi.mathnet.ru/eng/tm4007}
{E. M. Chirka, ``Equilibrium Measures on a Compact Riemann Surface'', \textit{Proc. Steklov Inst. Math.}, \textbf{306} (2019), 296--334.}

\bibitem{Chi19b} 
\href{http://mi.mathnet.ru/eng/mz12317}
{E. M. Chirka, ``Meromorphic Interpolation on a Compact Riemann Surface'', \textit{Math. Notes}, \textbf{106}:1 (2019), 156--159.}

\bibitem{Chi20} 
\href{http://mi.mathnet.ru/eng/book1806}
{E. M. Chirka, Capacities on a Compact Riemann Surface, in \textit{Analysis and mathematical physics}, Collected papers. On the occasion of the 70th birthday of Professor Armen Glebovich Sergeev, Tr. Mat. Inst. Steklova, \textbf{311}, ed. S. Yu. Nemirovski, A. V. Komlov, Steklov Math. Inst., Moscow, 2020, 281 pp.}

\bibitem{GoRa87} 
\href{http://mi.mathnet.ru/eng/msb2759}
{A. A. Gonchar, E. A. Rakhmanov, ``Equilibrium distributions and degree of rational
approximation of analytic functions'', \textit{Math. USSR-Sb.}, \textbf{62}:2 (1989), 305--348.}

\bibitem{GoRaSu11} 
\href{http://mi.mathnet.ru/eng/umn9452}
{A. A. Gonchar, E. A. Rakhmanov, S. P. Suetin,
``Pad\'e-Chebyshev approximants of multivalued analytic functions,
variation of equilibrium energy, and the S-property of stationary compact sets'',
\textit{Russian Math. Surveys}, \textbf{66}:6 (2011), 1015--1048.}

\bibitem{IkSu20} 
\href{http://mi.mathnet.ru/eng/tm4080}
{N. R. Ikonomov, S. P. Suetin, ``Scalar Equilibrium Problem and
the Limit Distribution of Zeros of Hermite--Pad\'{e} Polynomials of Type II'',
\textit{Modern problems of mathematical and theoretical physics}, Proc. Steklov Inst. Math., \textbf{309}, 2020, 159--182.}

\bibitem{KoPaSuCh17} 
\href{http://mi.mathnet.ru/eng/umn9786}
{A. V. Komlov, R. V. Palvelev, S. P. Suetin,
E. M. Chirka, ``Hermite--Pad\'e approximants for meromorphic functions on a compact
Riemann surface'', \textit{Russian Math. Surveys}, \textbf{72}:4 (2017), 671--706.}

\bibitem{Kom20} 
\href{https://www.springer.com/series/13332}
{A. V. Komlov, ``Polynomial Hermite--Pad\'{e} $m$-system and reconstruction of the values
of algebraic functions'', Submitted, 2020 (Research Perspectives CRM Barcelona, IRP
SAFAIS 2019), Trends in Mathematics, 2020, ISSN: 2509-7407.}

\bibitem{LoLo18} 
\href{https://doi.org/10.1016/j.jat.2017.10.002}
{A. L\'{o}pez-Garcia and G. L\'{o}pez Lagomasino, ``Nikishin systems on star-like sets: Ratio
asymptotics of the associated multiple orthogonal polynomials'',
\textit{J. Approx. Theory}, \textbf{225} (2018), 1--40.}

\bibitem{Lop20} 
\href{http://www.mathnet.ru/php/seminars.phtml?&presentid=29247&option_lang=eng}
{E. Lopatin, On a generalization of the scalar approach to the problem of limit distribution of zeros of Hermite--Pad\'{e} polynomials for a pair of functions forming the Nikishin system, Seminar on Complex Analysis (Gonchar Seminar), December 21, 2020.}

\bibitem{MaRaSu16} 
\href{http://mi.mathnet.ru/eng/conm2}
{Andrei Mart{\'\i}nez-Finkelshtein, Evguenii A. Rakhmanov, Sergey P. Suetin,
``Asymptotics of type I Hermite-Pad\'{e} polynomials for semiclassical functions'',
Modern trends in constructive function theory, \textit{Contemp. Math.}, \textbf{661}, 2016, 199--228.}

\bibitem{NuSi77} 
\href{http://dx.doi.org/10.1016/0021-9045(77)90117-4}
{J. Nuttall, R.S. Singh, ``Orthogonal polynomials and Pad\'{e} approximants associated
with a system of arcs'', \textit{J. Approx. Theory}, \textbf{21} (1977), 1--42.}

\bibitem{Nut84} 
\href{http://www.ams.org/mathscinet-getitem?mr=0769985}
{J. Nuttall, ``Asymptotics of diagonal Hermite--Pad\'e polynomials'', \textit{J. Approx.Theory},
\textbf{42} (1984), 299--386.}

\bibitem{Rak12} 
\href{https://doi.org/10.1090/conm/578/11484}
{E. A. Rakhmanov, ``Orthogonal polynomials and $S$-curves'', \textit{Recent advances in
orthogonal polynomials}, special functions and their applications, 11th International
Symposium, (August 29-September 2, 2011 Universidad Carlos III de Madrid Leganes,
Spain), Contemp. Math., \textbf{578}, eds. . J. Arvesu and G. Lopez Lagomasino, Amer.
Math. Soc., Providence, RI, 2012, 195--239.}

\bibitem{RaSu13} 
\href{http://mi.mathnet.ru/eng/msb8168}
{E. A. Rakhmanov, S. P. Suetin, ``The distribution of the zeros of the
Hermite-Pad\'e polynomials for a pair of functions forming a Nikishin system'',
\textit{Sb. Math.}, \textbf{204}:9 (2013), 1347--1390.}

\bibitem{Rak18} 
\href{http://mi.mathnet.ru/eng/umn9832}
{E. A. Rakhmanov, ``Zero distribution for Angelesco Hermite-Pad\'{e} polynomials'',
\textit{Russian Math. Surveys}, \textbf{73}:3 (2018), 457--518.}

\bibitem{ScSp54} 
\href{http://www.ams.org/mathscinet-getitem?mr=0065652}
{Menahem Schiffer, Donald C. Spencer, Functionals of finite
Riemann surfaces, Princeton University Press, Princeton, N. J., 1954, x+451 pp.}

\bibitem{Sor20} 
\href{http://mi.mathnet.ru/eng/msb8634}
{V. N. Sorokin, ``Hermite-Pad\'{e} approximations of Weyl function and
its derivative for discrete measures'', \textit{Mat. Sb.}, \textbf{211}:10 (2020).}

\bibitem{Sta88} 
\href{http://www.ams.org/mathscinet-getitem?mr=1005350}
{H. Stahl, ``Asymptotics of Hermite--Pad\'e polynomials and related convergence results.
A summary of results'', \textit{Nonlinear numerical methods and rational approximation}
(Wilrijk, 1987), Math. Appl., \textbf{43}, Reidel, Dordrecht, 1988, 23--53; also the fulltext
preprint version is avaible, 79 pp.}

\bibitem{Sta97b} 
\href{http://www.ams.org/mathscinet-getitem?mr=1484040}
{H. Stahl, ``The convergence of Pad\'{e} approximants to functions with branch points'',
\textit{J. Approx. Theory}, \textbf{91}:2 (1997), 139--204.}

\bibitem{Sta12} 
Herbert R. Stahl, \textit{Sets of Minimal Capacity and Extremal Domains}, 112 pp., arXiv: \href{http://arxiv.org/abs/1205.3811}{1205.3811}.

\bibitem{Sue00} 
\href{http://mi.mathnet.ru/eng/msb508}
{S. P. Suetin, ``Uniform convergence of Pad\'{e} diagonal approximants for hyperelliptic functions'',
\textit{Sb. Math.}, \textbf{191}:9 (2000), 1339--1373.}

\bibitem{Sue18} 
\href{http://mi.mathnet.ru/eng/mzm12181}
{S. P. Suetin, ``On an Example of the Nikishin System'',
\textit{Math. Notes}, \textbf{104}:6 (2018), 905--914.}

\bibitem{Sue18b} 
\href{http://mi.mathnet.ru/eng/tm3908}
{S. P. Suetin, ``On a new approach to the
problem of distribution of zeros of Hermite-Pad\'{e} polynomials for a Nikishin system'',
\textit{Proc. Steklov Inst. Math.}, \textbf{301} (2018), 245--261.}

\bibitem{Sue18d} 
Sergey P. Suetin, Hermite-Pad\'{e} polynomials and analytic continuation: new approach
and some results, 2018, 45 pp., arXiv: \href{http://arxiv.org/abs/1806.08735}{1806.08735}.

\bibitem{Sue19} 
\href{http://mi.mathnet.ru/eng/umn9884}
{S. P. Suetin, ``Existence of a three-sheeted Nutall surface for a certain class
of infinite-valued analytic functions'', \textit{Russian Math. Surveys}, \textbf{74}:2 (2019), 363--365.}

\bibitem{Sue19b} 
\href{http://mi.mathnet.ru/eng/mz12451}
{S. P. Suetin, ``Equivalence of a Scalar and a Vector Equilibrium Problem for
a Pair of Functions Forming a Nikishin System'', \textit{Math. Notes}, \textbf{106}:6 (2019), 971--980.}

\bibitem{Sue20} 
\href{http://mi.mathnet.ru/eng/umn9954}
{S. P. Suetin, ``Hermite--Pad\'{e} polynomials and Shafer quadratic approximations for
multivalued analytic functions'', \textit{Russian Math. Surveys}, \textbf{75}:4 (2020), 788--790.}

\end{thebibliography}
\end{document}